\documentclass[12pt,reqno,a4paper]{amsart}
\usepackage{color}
\usepackage{amssymb,amsmath}
\usepackage[latin1]{inputenc}
\usepackage[active]{srcltx}
\usepackage{graphicx}
\usepackage{exscale,relsize}
\usepackage{textgreek}
\usepackage{epsfig,graphics}
\usepackage{psfrag}
\usepackage{caption}
\usepackage{subcaption}
\usepackage[toc,page]{appendix}
\usepackage{bigints}

\usepackage[margin=1.4in]{geometry}

\def\N{\mathbb{N}}
\def\R{\mathbb{R}}
\def\m1{{I\!\!M}}

\renewcommand{\to}{\rightarrow}
\newcommand{\pa}{\partial}

\newcommand{\ino}{\int_{\Omega}}

\newcommand{\ainf}{\mbox{as\;}\;n\to+\infty}

\newcommand{\rife}[1]{(\ref{#1})}
\newcommand{\ov}[1]{\overline{#1}}
\newcommand{\un}[1]{\underline{#1}}

\newcommand{\sscp}{\scriptscriptstyle}

\renewcommand{\dfrac}{\displaystyle\frac}
\newcommand{\finedim}{\hspace{\fill}$\square$}
\newcommand{\intbar}{\mathop{\int\makebox(-15.5,0){\rule[6pt]{.7em}{0.3pt}}\kern-6pt}\nolimits}

\newcommand{\ii}{\infty}

\newcommand{\eps}{\varepsilon}
\newcommand{\dt}{\delta}

\newcommand{\sg}{\sigma}

\newcommand{\om}{\Omega}
\newcommand{\lm}{\lambda}

\newcommand{\vm}{v_{\sscp \mu}}

\renewcommand{\rho}{\mbox{\Large \textrho}}

\newcommand{\um}{u_{\sscp \mu}}

\newcommand{\vxi}{\xi}

\newtheorem{theorem}{Theorem}[section]
\newtheorem{proposition}[theorem]{Proposition}
\newtheorem{lemma}[theorem]{Lemma}
\newtheorem{corollary}[theorem]{Corollary}
\newtheorem{remark}[theorem]{Remark}
\newtheorem{definition}[theorem]{Definition}
\newcommand{\brm}{\begin{remark}\rm}
\newcommand{\erm}{\end{remark}}
\newcommand{\bdf}{\begin{definition}\rm}
\newcommand{\edf}{\end{definition}}
\newcommand{\bte}{\begin{theorem}}
\newcommand{\ete}{\end{theorem}}
\newcommand{\bpr}{\begin{proposition}}
\newcommand{\epr}{\end{proposition}}
\newcommand{\ble}{\begin{lemma}}
\newcommand{\ele}{\end{lemma}}
\newcommand{\bco}{\begin{corollary}}
\newcommand{\eco}{\end{corollary}}
\newcommand{\beq}{\begin{equation}}
\newcommand{\eeq}{\end{equation}}
\newcommand{\bdm}{\begin{displaymath}}
\newcommand{\edm}{\end{displaymath}}

\newcommand{\graf}[1]{\left\{\begin{array}{ll}#1\end{array}\right.}

\def\sideremark#1{\ifvmode\leavevmode\fi\vadjust{\vbox to0pt{\vss
 \hbox to 0pt{\hskip\hsize\hskip1em \vbox{\hsize2.1cm\tiny\raggedright\pretolerance10000 \noindent #1\hfill}\hss}\vbox to15pt{\vfil}\vss}}}

\begin{document}
\numberwithin{equation}{section}
\parindent=0pt
\hfuzz=2pt
\frenchspacing

\title[Generic properties of the Rabinowitz continuum]{
Generic properties of the Rabinowitz continuum}

\author[D. Bartolucci]{Daniele Bartolucci$^{(\dag)}$}
\address{Daniele Bartolucci, Department of Mathematics, University of Rome ``Tor Vergata", Via della ricerca scientifica 1, 00133 Roma, Italy}
\email{bartoluc@mat.uniroma2.it}

\author[Y. Hu]{Yeyao Hu}
\address{Yeyao Hu, School of Mathematics and Statistics, HNP-LAMA, Central South University,
Changsha, Hunan 410083, P.R. China}
\email{huyeyao@csu.edu.cn}

\author[A. Jevnikar]{Aleks Jevnikar}
\address{Aleks Jevnikar, Department of Mathematics, Computer Science and Physics, University of Udine, Via delle Scienze 206, 33100 Udine, Italy}
\email{aleks.jevnikar@uniud.it}

\author[W. Yang]{Wen Yang}
\address{Wen Yang, Wuhan Institute of Physics and Mathematics, Chinese Academy of Sciences, P.O. Box 71010,
Wuhan 430071, P.R. China}
\address{Innovation Academy for Precision Measurement Science and Technology, Chinese Academy of
Sciences, Wuhan 430071, P.R. China}
\email{wyang@wipm.ac.cn}


\thanks{$^{(\dag)}$Research partially supported by:
Beyond Borders project 2019 (sponsored by Univ. of Rome "Tor Vergata") "{\em Variational Approaches to PDE's}",
MIUR Excellence Department Project awarded to the Department of Mathematics, Univ. of Rome Tor Vergata, CUP E83C18000100006.}

\begin{abstract}
In this paper we prove that generically, in the sense of domain variations, the unbounded Rabinowitz continuum of solutions to a nonlinear eigenvalue problem is a simple analytic curve. The global bifurcation diagram resembles the classic model case of the Gel'fand problem in dimension two.
\end{abstract}
\maketitle



\

\setcounter{section}{0}
\section{\bf Introduction}
\setcounter{equation}{0}

Let $\om\subset \R^N$ be an open and bounded domain of class $C^4$,
we are concerned with generic properties of the Rabinowitz (\cite{Rab}) unbounded continuum of
$C^{2,r}_0(\ov{\om})$-solutions of
\beq\label{pm}
\graf{-\Delta v=\mu f(v) \quad \mbox{in}\;\;\om\\ v=0 \quad \mbox{on}\;\;\pa \om}
\eeq
with $\mu\geq 0$ and $f$ satisfying:\\

${\bf (H1)}$ $f:(a,+\ii)\to (0,+\ii)$ of class $C^2$ for some $a<0$,
$f^{'}(t)>0, f^{''}(t)>0$, $\forall\, t\in (a,+\ii)$.\\

In particular by the maximum principle we have $v>0$ in $\om$.
It follows from \cite{Rab} that there exists a closed
(in the $[0,+\ii)\times C^{2,r}_0(\ov{\om})$-topology) connected and unbounded
set of solutions $(\mu,\vm)$ of \rife{pm}, which we denote by $\mathcal{R}_{\ii}$, that contains
the unique solution for $\mu=0$, which is $(\mu,\vm)=(0,0)$. Of course it is not true in general
that $\mathcal{R}_{\ii}$ is a simple curve with no
bifurcation points, see for example \cite{N,NS,PW}. If $f$ is real analytic,
then $\mathcal{R}_{\ii}$ is also a path-connected set (\cite{Dan3}). On the other side, much more is known for certain classes
of nonlinearities in the radial case (see \cite{JL} and in particular \cite{Dan1,Kor,Kor1} for a review)
or, limited to $\om\subset \R^2$, for symmetric and convex geometries (see \cite{HK}). Actually,
in these cases in particular $\mathcal{R}_{\ii}$ is a $1$-dimensional connected manifold in
$[0,+\ii)\times C^{2,r}_0(\ov{\om})$ whose boundary is $(0,0)$.
See also \cite{BJ} for a more detailed description of the \un{qualitative} behavior of $\mathcal{R}_{\ii}$
for $f(t)=e^t$ and $N=2$.\\
Our aim here is to prove, under suitable regularity and growth condition on $f$, that for ``almost any domain" in a suitably defined sense,
$\mathcal{R}_{\ii}$ is indeed a $1$-dimensional connected manifold in
$[0,+\ii)\times C^{2,r}_0(\ov{\om})$ whose boundary is $(0,0)$. Although
this result seems to be well known, we could not find a statement of this sort in
literature.\\

As far as we are concerned with generic properties with respect to domain variations,
many results are by now classical, see for example \cite{ST}, \cite{H} and references quoted therein.
Among many other applications which we cannot discuss here,
it follows from \cite{ST}, \cite{H} that, for a fixed $\mu$,
then for ``almost any domain" any solution of \rife{pm} is non degenerate.
This is obviously false in general if $\mu$ is not fixed (\cite{CrRab}).
Actually these sort of results hold for much more general
semilinear elliptic PDE's and are used (\cite{ST}, \cite{H}) to infer that
the number of solutions of certain equations is, generically with respect to domain variations,
either finite or at most countable.\\
Generic simpleness of eigenvalues and/or non-degeneracy
properties with respect to variations of $\mu$ and or coefficients of the equations
are also well known, starting with \cite{Mich,Uhl,ST} (see also \cite{Rynne}) as later
improved in the real-analytic framework in \cite{Dan}.\\
None of these results seems to fit our need as we seek a property which guarantees, in a
generic sense with respect to domain variations, that $\mathcal{R}_{\ii}$
is a simple curve, that is in particular, without self intersections or bifurcation points.
It seems that this kind of result neither follows from the highly sophisticated
global analytic bifurcation theory (\cite{bdt,but}), which, in its full generality,
does not guarantee the global injectivity of the branch.\\

While it is not true that any solution $(\mu,\um)$ of \rife{pm} is, generically with respect
domain variations, non degenerate (\cite{CrRab}),
on the other side we can prove that, in a generic sense, either $(\mu,\um)$  is non degenerate or the
classical Crandall-Rabinowitz (\cite{CrRab}) bending condition is satisfied (see Proposition \ref{pr3.1} below).

The argument works essentially as in \cite{ST}, but we find that
some major simplifications are obtained by taking the point of view of \cite{H}, where
one can find a beautifully refined and simplified transversality theory with respect to domain variations.\\

For $\om_0$ a bounded domain of class $C^4$ (see section \ref{sec3} for details)
we denote by $\mbox{\rm Diff}^{\,4}(\om_0)$ the set of diffeomorphisms
$h:\ov{\om_0}\to\ov{\om}$ of class $C^4$. We recall that a subset of a metric space is said to be:\\
- nowhere dense, if its closure has empty interior; \\
- meager (or of first Baire category), if it is the union of countably many nowhere dense sets.\\

Once more, it is likely that this result is known to experts in the field, still we could not find a
statement of this sort in literature.\\
Here $L_{\mu}$ is the linearized operator relative to \rife{pm} (see section \ref{sec2}).
Then we have

\bte\label{thm1.generic} Let $f:(a,+\ii)\to (0,+\ii)$ be of class $C^2$ for some $a<0$.
For any $\om_0\subset \R^N$ of class $C^{4}$ there exists a meager set
$\mathcal{F}\subset \mbox{\rm Diff}^{\,4}(\om_0)$, depending on $f,N,\om_0$,
such that if $h\in \mbox{\rm Diff}^{\,4}(\om_0)\setminus\mathcal{F}$ then, for any
solution $(\mu,\vm)$ of \mbox{\rm \rife{pm}} on $\om:=h(\om_0)$ with $\mu>0$, it holds: either\\
$(a)$   Ker$(L_{\mu})=\emptyset$, or\\
$(b)$  Ker$(L_{\mu})=\mbox{\rm span}\{\phi\}$ is one dimensional and
$\ino f(\vm)\phi\neq 0$.
\ete

In particular we deduce the following result about the Rabinowitz (\cite{Rab}) unbounded
continuum of solutions of \rife{pm} for $f$ real analytic, which satisfies {\bf (H1)} and,\\

{\bf (H2)} for any $\dt>0$ there exists $C_{\dt}>0$ (depending also by $f, N, \om$) such that
$\vm\leq C_{\dt}$ for any solution of \rife{pm} with $\mu\geq \dt$.\\

It is well known that, for $\om$ a bounded domain of class $C^4$, ${\bf (H2)}$
is satisfied under suitable growth assumptions on $f$, as for example those in \cite{DLN},
(here $F(t)=\int\limits_{0}^{t}f(s)ds$),
$$
\graf{\lim\limits_{t\to +\ii}\frac{f(t)}{t}=+\ii,\quad \lim\limits_{t\to +\ii}
\frac{f(t)}{t^{\beta}}=0,\quad
\beta=\frac{N+2}{N-2}\,\mbox{ if } N\geq 3,\;\beta<+\ii\,\mbox{ if } N=2,\\
\limsup\limits_{t\to +\ii} \frac{tf(t)-\theta F(t)}{t^2 f^{\frac{2}{N}}(t)}\leq 0,
\quad \mbox{ for some }\theta \in [0,\frac{2N}{N-2}),\\
\mbox{ and if }N\geq 3,\; f(t)t^{-\frac{N+2}{N-2}}\;\mbox{ is non increasing in }(0,+\ii)}
$$

Clearly the model nonlinearities $f(t)=(1+t)^p$, $1<p<\frac{N+2}{N-2}$, $N\geq 3$, $p>1$, $N=2$, fit 
these assumptions. However there are many other cases where ${\bf (H2)}$ is satisfied, as for example $f(t)=e^t$, $N=2$ (\cite{NS90}). Then we have,

\bte\label{Rab}
Let $f$ be real analytic and satisfying {\bf (H1)-(H2)}, $\om_0\subset \R^N$ of class $C^4$,
$\mathcal{F}\subset \mbox{\rm Diff}^{\,4}(\om_0)$ as defined by Theorem \ref{thm1.generic} and
pick $h\in \mbox{\rm Diff}^{\,4}(\om_0)\setminus\mathcal{F}$. If $\om=h(\om_0)$, then the Rabinowitz
unbounded continuum $\mathcal{R}_\ii$ of solutions of {\rm \rife{pm}} is a 1-dimensional
real analytic manifold with boundary $(\mu(0),v(0))=(0,0)$. In particular
$$
\mathcal{R}_{\ii}=\Bigr\{[0,\ii)\ni s \mapsto (\mu(s),v(s))\in [0,+\ii)\times C^{2,r}_0(\ov{\om}\,)\Bigr\},
$$
is a continuous simple curve without bifurcation points
where $(\mu(0),v(0))=(0,0)$ and $\mu(s)\to 0^+$ and $\|v(s)\|_\ii\to +\ii$ as $s\to \ii$.
\ete

\brm \emph{Theorem \ref{thm1.generic} can be generalized to the case of uniformly elliptic operators such as $Lu:=\textup{div}(A(\nabla u))+ \vec{b}\cdot\nabla u + c u$, where $A=(a_{ij})$, $a_{ij}(x)$, $b_j(x)$ and $c(x)$ are smooth up to the boundary. Therefore, the generic bending result (Theorem \ref{Rab}) also follows if one replace the Laplace operator in \eqref{pm} by uniformly elliptic operators.}
\erm

\bigskip
\bigskip

{\bf Acknowledgments.} The first author would like to express His warmest thanks to
N. Dancer for pointing out that a combined use of Lemmas 9 and 10 in \cite{Dan3} shows
that $\mathcal{R}_\ii$ is path-connected and to B. Buffoni for
very useful discussions about analytic global bifurcation theory.

\bigskip
\bigskip

\section{\bf Well know results}\label{sec2}
Let $X=\R\times C^{2,r}_0(\ov{\om})$, we introduce the map,
\beq\label{eF}
F: X \to  C^{r}(\ov{\om}),\quad F(\mu,v):=-\Delta v -\mu f(v),
\eeq
and its differential with respect to $(\mu,v)$, that is the linear operator,
$$
D_{\mu,v}F(\mu,v):X \to  C^{r}(\ov{\om}\,),
$$
which acts as follows,
$$
D_{\mu,v}F(\mu,v)[\dot{\mu},\dot{v}]=
D_v F(\mu,v)[\dot{v}]+d_{\mu}F(\mu,v) [\dot{\mu}],
$$
where we have introduced the differential operators,
$$
D_v F(\mu,v)[\dot{v}]=
-\Delta \dot{v} -\mu f^{'}(v) \dot{v}, \quad \dot{v} \in C^{2,r}_0(\ov{\om}\,),
$$

$$
d_{\mu}F(\mu,v) [\dot{\mu}]=-f(v)\dot{\mu}, \quad \dot{\mu} \in \R.
$$

For a fixed solution $(\mu,\vm)$ the eigenvalues of $L_\mu:=D_v F(\mu,\vm)$ form an increasing sequence and
are denoted by $\sg_k$, $k\in \N$, which satisfy
$$
L_{\mu}\phi=\sg_k\phi,\phi\in C^{2,r}_0(\ov{\om}).
$$

\bigskip

By the Fredholm alternative,
the implicit function theorem applies around any solution of \rife{pm} as follows:
\ble\label{lem1.1}
Let $(\mu_0,v_0)$ be a solution of {\rm \rife{pm}} with $\mu=\mu_0\geq 0$.\\
If $0$ is not an eigenvalue of $L_{\sscp \mu_0}$, then:\\
$(i)$ $L_{\sscp \mu_0}$ is an isomorphism;\\
$(ii)$ There exists an open neighborhood $J\subset \R$ of $\mu_0$ and
$\mathcal{B}\subset C^{2,r}_{0}(\ov{\om}\,)$ of $v_0$, such
that the set of solutions of
{\rm \rife{pm}} in $J\times \mathcal{B}$ is a curve of class $C^{2}$, $J\ni\mu \mapsto \vm \in \mathcal{B}$.
\ele

Next we state the well known bending result of \cite{CrRab} for solutions of $\rife{pm}$ just
with an additional observation about the case where $f$ is real analytic in $(a,+\ii)$ for some $a<0$.
The conclusions deduced in this particular case are straightforward consequences of general
and well known facts of analytic bifurcation theory (\cite{but}).

\bpr{{\rm (\cite{CrRab})}}\label{pr3.1}
Let $(\mu,\vm)$ be a solution of {\rm \rife{pm}} with $\mu>0$ and suppose that the  $k$-th
eigenvalue of $L_{\mu}$ satisfies $\sg_k=0$ and is simple, that is, it admits only one eigenfunction,
$\phi_k\in C^{2,r}_0(\ov{\om}\,)$. If
$$
\ino f(\vm)\phi_k\neq 0,
$$
then there exists $\eps>0$, an open neighborhood $\mathcal{U}$ of $(\mu,\vm)$ in $X$ and
a curve $(-\eps,\eps)\ni s \mapsto (\mu(s),v(s))$ of class $C^2$ such that
$(\mu(0),v(0))=(\mu,\vm)$ and the set of solutions of {\rm \rife{pm}} in $\mathcal{U}$
has the form $(\mu(s),v(s))$ with,
$v(s)=\vm+s\phi_k+\vxi(s)$, and
$$
\ino f(v(s))\vxi(s)\phi_k=0,\quad s\in (-\eps,\eps).
$$
Moreover it holds,
\beq\label{2907.1}
\vxi(0)\equiv 0\equiv \vxi^{'}(0),\quad\mu^{'}(0)=0,
\eeq
and there exists a continuous curve $(\sg(s),\phi(s))$, such that
$L_{\mu(s)}\phi(s)=\sg(s)\phi(s)$, $s\in (-\eps,\eps)$, $\phi(0)=\phi_k$, $\sg(0)=\sg_k$ and
$$
\sg(s)\ino f^{'}(v(s))\phi(s)v(s) \quad \mbox{and} \quad \mu^{'}(s)\ino f(v(s))\phi(s)
$$
have the same zeroes and,
whenever $\mu^{'}(s)\neq 0$, the same sign. In particular
$$
\dfrac{\sg(s)}{\mu^{'}(s)}=\dfrac{\ino f(\vm)\phi_k+\mbox{\rm o}(1)}
{\ino \phi^2_k+\mbox{\rm o}(1)},\mbox{ as }s\to 0.
$$

If $f$ is real analytic in $(a,+\ii)$ for some $a<0$ then $\mu(s),v(s),\sg(s),\phi(s)$ are real
analytic functions of $s\in (-\eps,\eps)$ and in particular either $\mu(s)$
is constant in $(-\eps,\eps)$ or $\mu^{'}(s)\neq 0$, $\sg(s)\neq 0$ in $(-\eps,\eps)\setminus\{0\}$
and $\sg(s)$ is simple in $(-\eps,\eps)$.
\epr

\bigskip
\bigskip

\section{\bf Generic properties of the Rabinowitz continuum}\label{sec3}
In this section we prove Theorem \ref{thm1.generic} and \ref{Rab}.\\

Let us recall few definitions and set some notations first.
\bdf\label{C1} {\it
A domain $\om$ is of class $C^{k}(C^{k,r})$, $k\geq 1$, if for each $x_0\in \pa\om$ there exists a ball
$B=B_r(x_0)$
and a one to one map $\Theta: B\mapsto U\subset \R^N$ such that $\Theta\in C^{k}(B)(C^{k,r}(B)),
\Theta^{-1}\in C^{k}(U)(C^{k,r}(U))$ and the following holds:
$$
\Theta(\om\cap B)\subset \R^N_+\quad \mbox{ and }\quad \Theta(\om\cap B)\subset \pa\R^N_+.
$$

It is well known {\rm(}see for example \cite{H}{\rm)} that this is equivalent to say that there exists
$r>0$ and $M>0$ such that, given any ball
$B_r(x_0)$, $x_0\in \R^2$ then, after suitable rotation and translations, it holds:
$$
\om\cap B=\{(x_1,x_2)\,:\,x_2<f(x_1)\}\cap B
$$
and
$$
\pa\om\cap B=\{(x_1,x_2)\,:\,x_2=f(x_1)\}\cap B,
$$
for some $f\in C^{k}(\R)(C^{k,r}(\R))$ whose norm is not larger than $M.$}
\edf

\medskip

\bdf {\it
Let $\om\subset \R^N$ be an open and bounded domain of class $C^m$, $m\geq 1$.
$C^m(\ov{\om}\,;\R^N)$ is the Banach space of continuous and $m$-times differentiable maps on $\om$,
whose derivatives of order $j=0,1,\cdots,m$
extend continuously on $\ov{\om}$. $\mbox{\rm Diff}^{\,m}(\om)\subset C^m(\ov{\om}\,;\R^N)$ is the open subset of
$C^m(\ov{\om}\,;\R^N)$
whose elements are $C^m$ imbeddings on $\ov{\om}$, that is, of maps $h:\ov{\om}\mapsto \R^N$ which are diffeomorphisms of class $C^m$
on their images ${h(\ov{\om}\,)}$.}
\edf

We recall that if $X,Z$ are Banach spaces and $T:X\to Z$ is linear and continuous, then T is Fredholm (semi-Fredholm)
if $R(T)$ (the range of $T$) is closed and both dim(Ker$(T)$) and codim$(R(T))$ are finite.
If $T$ is Fredholm, then the index of $T$ is
$$
\mbox{\rm ind}(T)=\mbox{\rm dim(Ker$(T)$)}-\mbox{\rm codim}(R(T)).
$$
We refer to \cite{Ka} for further details about Fredholm operators. Given a Banach space
$X$ and $x\in X$, we will denote by $T_x X$ the tangent space at $x$.

\bdf{\it
Let $X,Z$ be Banach spaces, $A\subset X$ an open set and $F:A\to Z$ a $C^1$ map. Suppose that for any $x \in  A$ the
Fr\'echet derivative $D_x F(x):T_xX\to T_\eta Z$ is a Fredholm operator. A point $x\in A$ is a \un{regular point}
if $D_x F(x)$ is surjective, is a \un{singular point} otherwise. The image of a singular point $\eta=F(x)\in Z$ is a \un{singular value}.
The complement of the set of singular values in $Z$ is the set of \un{regular values}.}
\edf

The following Theorem is a particular case of a more general transversality result proved in \cite{H}, see also \cite{Sm}.

\bte[{\rm \cite{H}}]\label{enry}
Let ${X},\mathcal{H},Z$ be separable Banach spaces,
$\mathcal{A}\subseteq X\times \mathcal{H}$ an open set,
${\Phi}:\mathcal{A}\to Z$ a map of class $C^k$ and
$\eta\in Z$.\\ Suppose that for each $(x,h)\in {\Phi}^{-1}(\eta)$ it holds:
\begin{align}\nonumber
&(i) D_{x}{\Phi}(x,h):T_x X\to T_\eta Z\;\mbox{is a Fredholm operator with index $< k$};\\
&(ii) D {\Phi}(x,h)=(D_{x}{\Phi}(x,h),D_{h}{\Phi}(x,h)):
T_x X\times T_h\mathcal{H} \to T_\eta Z\;\mbox{is surjective}.\nonumber
\end{align}
Let $A_h=\{x\,:\,(x,h)\in \mathcal{A}\}$ and
$$
\mathcal{H}_{\rm crit}=\{h\,:\,\eta \mbox{ is a singular value of }{\Phi}(\,\cdot\,, h):A_h\to Z\}.
$$
Then $\mathcal{H}_{\rm crit}$ is meager in $\mathcal{H}$.
\ete

\bigskip
\bigskip

We are ready to present the proof of Theorem \ref{thm1.generic}.\\
{\it Proof of Theorem \ref{thm1.generic}}\\
Let $\om_0$ as in the statement and let us define
$$
X_{\sscp \om_0}=\R \times C^{2,r}_0(\ov{\om_0}\,).
$$
We define the maps,
$$
F_{\om_0}:X_{\sscp \om_0} \to C^r(\ov{\om_0}\,  ),
\quad  F_{\om_0}(\mu,v)=\Delta v+ \mu f(v).
$$

Next, for fixed $h\in \mbox{\rm Diff}^{\, 4}(\om_0)$ and
$v\in C^{2,r}_0(\ov{h(\om_0)}\,)$,
we define the pull back,
$$
h^*(v)(x)=v(h(x)),\;x\in \ov{\om_0}.
$$
Clearly $h^*$ is an isomorphism of $C^{2,r}_0(\ov{h(\om_0)}\,)$ onto $C^{2,r}_0(\ov{\om_0}\,)$
with inverse $h^{*-1}=(h^{-1})^*$. For any such $h$, it is well defined the map
$$
F_{h(\om_0)}:X_{h(\om_0)}\to C^{r}(\ov{h(\om_0)}\,)
$$
and then we can set,
$$
h^*F_{h(\om_0)}h^{*-1}:X_{\om_0}\times \mbox{\rm Diff}^{\, 4}(\om_0)\to C^{r}(\ov{\om_0}\,).
$$

Putting $\mathcal{H}=\mbox{\rm Diff}^{\, 4}(\om_0)$,
$\eta=0\in Z=C^{r}(\ov{\om_0}\,)$, we will apply Theorem \ref{enry} to the map
${\Phi}={\Phi}(\mu,v,h)$ defined as follows

$$
{\Phi}:\mathcal{A}\to \R\times C^{r}(\ov{\om_0}\,),\qquad
\mathcal{A}=X_{\sscp \om_0}\times \mathcal{H},
$$

$$
{\Phi}(\mu,v,h)=h^*F_{h(\om_0)}h^{*-1}(\mu,v)
$$

\bigskip
\bigskip
{\bf STEP 1:} Our aim is to show that the assumptions $(i)$ and $(ii)$ of Theorem~\ref{enry} hold.\\
As in \cite{H}, it is very useful for the discussion to denote by
$(\dot{\mu},\dot{v},\dot{h})\in \R\times C^{2,r}_0(\ov{\om_0}\,)\times C^{4}(\ov{\om_0}\,;\,\R^N)$ the
elements of the tangent space at points
$({\mu},v,{h})\in X_{\sscp \om_0} \times \mathcal{H}$.\\
First of all observe that for fixed $h\in \mbox{\rm Diff}^{\, 4}(\om_0)$,
the linearized  operator,
$$
D_{\mu,v}{\Phi}(\mu,v,h):\R\times C^{2,r}_0(\ov{\om_0}\,) \to C^{r}(\ov{\om_0}\,),
$$
acts as follows on $(\dot{\mu},\dot{v})\in \R\times C^{2,r}_0(\ov{\om_0}\,)$,
$$
D_{\mu,v}{\Phi}(\mu,v,h)[\dot{\mu},\dot{v}]=
h^*\left(\Delta \dot{v}^*+\mu f^{'}(v^*)\dot{v}^*+f(v^*)\dot{\mu}\right),
$$
where
$$
{v}^*=(h^{*})^{-1}{v},\quad \dot{v}^*=(h^{*})^{-1}\dot{v}.
$$
Since any diffeomorphism of class $C^{4}$ maps the Laplace operator to a uniformly elliptic operator
with $C^2$ coefficients, by standard elliptic estimates it is not difficult  to see that
$D_{\mu,v}{\Phi}(\mu,v,h)$ is a Fredholm operator
of index 1.\\
This fact proves $(i)$ whenever we can show that ${\Phi}\in C^{k}(\mathcal{A})$ for some
$k\geq 2$.  The regularity of ${\Phi}$ with respect to $h$
is the same as that of $F_{h(\om_0)}$ with respect to $v$,
see chapter 2 in \cite{H}. Therefore,
we have ${\Phi}\in C^{3}(\mathcal{A})$, as claimed.\\

Next we prove $(ii)$, that is, we show that $\eta=0$ is a regular value for the map
$(\mu,v,h)\to {\Phi}(\mu,v,h)$.
We argue by contradiction and
suppose that there exists a singular point
$(\ov{\mu},\ov{v},\ov{h}\,)$ of ${\Phi}$ such that
${\Phi}(\ov{\mu},\ov{v},\ov{h}\,)=0$.\\
First of all, let us define $\om=\ov{h}(\om_0)$,
$\ov{u}=(\ov{h}^*)^{-1}\ov{v}\in C^{2,r}_0(\ov{\om}\,)$ and
$\widehat{\Phi}(\mu,u,h)$
on $X_{\sscp \om}\times \mbox{\rm Diff}^{\, 4}(\om)$ as follows,
$$
\widehat{\Phi}(\mu,u,h)=
h^*F_{\om}h^{*-1}(\mu,u)
$$
where,
$$
F_{\om}:X_{\sscp \om}\to C^r(\ov{\om}),
\quad  F_{\om}(\mu,u)=\Delta u+ \mu f(u),
$$
Let $i_\om\in  \mbox{\rm Diff}^{\, 4}(\om)$ be the identity map.
By construction, in these new coordinates the map
$\widehat{\Phi}(\mu,u,h)$ has a singular point $(\ov{\mu},\ov{u},i_\om)$ such that
$\widehat{\Phi}(\ov{\mu},\ov{u},i_\om)=0$, that is,
by assumption the derivative $D_{\mu,u,h}\widehat{\Phi}(\ov{\mu},\ov{u},i_\om)$
is not surjective. Putting
$$
\ov{f}=f(\ov{u}), \quad \ov{f}^{'}=f^{'}(\ov{u}),
$$
a subtle evaluation shows that
$D_{\mu,u,h}\widehat{\Phi}(\ov{\mu},\ov{u}, i_\om)$ acts on
$$
(\dot{\mu},\dot{u},\dot{h})\in \R\times C^{2,r}_0(\ov{\om})\times C^{4}(\ov{\om}\,;\R^N)
$$
as follows (see Theorem 2.2 in \cite{H}),
$$
D_{\mu,u,h}\widehat{\Phi}(\ov{\mu},\ov{u}, i_\om)[\dot{\mu},\dot{u},\dot{h}]=
$$
$$
\Delta \dot{u}+\ov{\mu}\ov{f}^{'} \dot{u} +\ov{f}\dot{\mu}
+
\dot{h}\cdot \nabla(\Delta \ov{u}+\ov{\mu}\ov{f})-
(\Delta +\ov{\mu}\ov{f}^{'})\dot{h}\cdot \nabla \ov{u}=
$$
\beq\label{DF1}
\left(\Delta +\ov{\mu}\ov{f}^{'}\right)\dot{u}-
\left(\Delta +\ov{\mu}\ov{f}^{'}\,\right)\dot{h}\cdot\nabla
\ov{u}+\ov{f}\dot{\mu},
\eeq
where we used the fact that
$\Delta \ov{u}+\ov{\mu}\ov{f}=\widehat{\Phi}(\ov{\mu},\ov{u}, i_\om)=0$.\\
At this point observe that, by the Fredholm property of the operator $\Delta +\ov{\mu}\ov{f}$ on $C^{2,r}_0(\ov{\om}\,)$,
we have that the subspace
$\left\{D_{\mu,u,h}\widehat{\Phi}(\ov{\mu},\ov{u}, i_\om)[(0,\dot{u},0)],
\,\dot{u}\in C^{2,r}_0(\ov{\om})\right\}$, is closed and has finite codimension. Next,
since $\ov{u}\in C^{2,r}_0(\ov{\om}\,)$ and $\pa\om$ is of class $C^{4}$, then
by standard elliptic regularity theory we find that $\ov{u}\in C^{3,r}_0(\ov{\om})$ and then
$\dot{h}\cdot \nabla \ov{u}\in C^{2,r}(\ov{\om})$. As a consequence we can prove
that  the subspace
$\left\{ D_{\mu,u,h}\widehat{\Phi}(\ov{\mu},\ov{u}, i_\om)[(0,0,\dot{h})],
\dot{h}\in C^{4}(\ov{\om}\,;\R^N)\right\}$ is closed with finite codimension as well.
Indeed, let us define
$K:C^{2,r}(\ov{\om}\,)\mapsto C^{2,r}(\ov{\om}\,)$ as the
linear operator which, to any $\phi\in C^{2,r}(\ov{\om}\,)$, associates the unique solution
$\phi_b=K[\phi]$ of $\Delta \phi_b=0$,
$\phi_b=\phi$ on $\pa\om$.
Clearly this is always well posed since $\om$ is of class $C^4$ and $\phi\in C^{2,r}(\ov{\om}\,)$. Then,
$$
\Delta \phi + \ov{\mu}\ov{f}^{'}\phi=g\in C^{r}(\ov{\om}\,),
$$
if and only if
$$
\phi \in C^{2,r}(\ov{\om}\,)\mbox{ and }\phi+T[\phi]=G[g]\in C^{2,r}(\ov{\om}\,),
$$
where $G[g]=\ino G(x,y)g(y)$ and $T:C^{2,r}(\ov{\om}\,)\mapsto C^{2,r}(\ov{\om}\,)$,
$T(\phi)=G[\ov{\mu}\ov{f}^{'}\phi]-K[\phi]$.
Since $\om$ is of class $C^{4}$, then by standard elliptic estimates, $T$
maps $C^{2,r}(\ov{\om}\,)$ into $C^{3,r}(\ov{\om}\,)$. Therefore $T$ is compact and
then we conclude by the Fredholm alternative that the range of
$(\Delta + \ov{\mu}\ov{f}^{'})(\dot{h}\cdot \nabla \ov{u})$,
$\dot{h}\cdot \nabla \ov{u}\in C^{2,r}(\ov{\om}\,)$,
is closed in $C^{r}(\ov{\om}\,)$ and has finite codimension.\\

At this point, we deduce from these two facts that there exists a non trivial
$\phi_{\perp}\in C^r(\ov{\om}\,)$ which is orthogonal to the image of
$D_{\mu,u,h}\widehat{\Phi}(\ov{\mu},\ov{u}, i_\om)$, that is,
\beq\label{1perp}
\ino \phi_{\perp}
\left(\left(\Delta +\ov{\mu}\ov{f}^{'}\right)\dot{u}-
\left(\Delta +\ov{\mu}\ov{f}^{'}\,\right)\dot{h}\cdot\nabla
\ov{u}+\ov{f}\dot{\mu}\right)=0,
\forall\,(\dot{\mu},\dot{u},\dot{h}).
\eeq
Putting $(\dot{\mu},\dot{u},\dot{h})=(\dot{\mu},0,0)$ in \rife{1perp}
we find $\ino \ov{f}\phi_{\perp}=0$, and then
if we choose  $\dot{h}=0$, we find that,
$$
\ino \phi_{\perp}\left(\Delta +\ov{\mu}\ov{f}^{'}\right)\dot{u}=0,\,\forall \, \dot{u}\in C^{2,r}_0(\ov{\om}\,),
$$
which shows that $\phi_{\perp}$ is a $C^{r}(\ov{\om}\,)$ distributional solution of
$\Delta  \phi_{\perp}+\ov{\mu}\ov{f}^{'}\phi_{\perp}=0$.
Therefore, by standard elliptic estimates (where we recall that $\pa\om$ is of class $C^{4}$),
$\phi_{\perp}$ is a $C_0^2(\ov{\om}\,)$ solution of
$\Delta  \phi_{\perp}+\ov{\mu}\ov{f}^{'}\phi_{\perp}=0$.
As a consequence we observe that \rife{1perp} is reduced to
$$
\ino \phi_{\perp}\left(\Delta +\ov{\mu}\ov{f}^{'}\right)\dot{h}\cdot \nabla\ov{u}=0,
\quad \forall\, \dot{h}\in C^{4}(\ov{\om}\,;\R^N),
$$
which allows us to deduce that,
$$
0=\ino \phi_{\perp}\left(\Delta +\ov{\mu}\ov{f}^{'}\right)\dot{h}\cdot \nabla\ov{u}=
$$
$$
\ino \phi_{\perp}\left(\Delta +\ov{\mu}\ov{f}^{'}\right)\dot{h}\cdot \nabla \ov{u}-
\ino \left(\Delta \phi_{\perp}+\ov{\mu}\ov{f}^{'}\phi_{\perp}\right)\dot{h}\cdot \nabla \ov{u}=
$$
$$
\ino \phi_{\perp}\Delta(\dot{h}\cdot \nabla \ov{u})-
\ino (\Delta\phi_{\perp})\dot{h}\cdot \nabla \ov{u}=
\int\limits_{\pa\om} \left(\phi_{\perp}\pa_{\nu}(\dot{h}\cdot \nabla \ov{u})-
\dot{h}\cdot \nabla \ov{u}(\pa_\nu\phi_{\perp})\right)=
$$
$$
-\int\limits_{\pa\om} (\pa_\nu\phi_{\perp})\dot{h}\cdot \nabla \ov{u}=
-\int\limits_{\pa\om} (\pa_\nu\phi_{\perp})(\pa_\nu \ov{u})\dot{h}\cdot \nu,\quad \forall\,
\dot{h}\in C^4(\ov{\om},\R^N).
$$

Therefore, since $\dot{h}$ is arbitrary, we conclude that,
$$
(\pa_\nu\phi_{\perp})(\pa_\nu \ov{u})\equiv0 \quad \mbox{ on }\pa\om.
$$
At this point we observe that since $\ov{f}>0$  on $\ov{\om}$ and $\ov{u}=0$ on $\pa\om$,
then, by the strong maximum principle,
we have $\ov{u}>0$ in $\om$. Since $\pa\om$ is of class $C^{4}$
we can apply the Hopf boundary Lemma and conclude that
$\pa_\nu \ov{u}<0$ on $\pa\om$. Therefore we conclude that necessarily
$\pa_\nu\phi_{\perp}\equiv 0$ on $\pa\om$, which is in contradiction
with the Hopf boundary Lemma. This contradiction shows that $(ii)$ holds
and then we can apply Theorem \ref{enry} and conclude that
there exists a meager set $\mathcal{F}\subset \mbox{Diff}^4(\om_0)$
such that if $h(\om_0)\notin \mathcal{F}$ then $\eta=0$ is a regular
value of ${\Phi}(\mu,v,h)$.\\

{\bf STEP 2:}
We have from STEP 1 that there exists a meager set $\mathcal{F}\subset \mbox{Diff}^4(\om_0)$
such that if $h\notin \mathcal{F}$ and
$\om:=h(\om_0)$,
then $\eta=0$ is a regular value of the map ${\Phi}(\cdot,\cdot,h)$.
As a consequence, for any $(\ov{\mu},\ov{v})$ which solves

$$
\Phi(\mu,v)=F_{\om}(\mu,v)=0
$$
and setting $\ov{f}=f(\ov{v})$, $\ov{f}^{'}(\ov{v})$ then the differential
$$
\ov{L}[\dot{\mu},\dot{v}]:=D_{\mu,v}\Phi(\ov{\mu},\ov{v})[\dot{\mu},\dot{v}] =
                        \Delta\dot{v}  +\ov{\mu} \ov{f}^{'}(\ov{v})\dot{v} + \ov{f}\dot{\mu}.
$$
is surjective.
On the other side, since $\ov{v}$ solves \rife{pm}, then the operator,
$$
\Delta\dot{v}  +\ov{\mu} \ov{f}^{'}(\ov{v})\dot{v}
$$
is just $L_{\ov{\mu}}$ for which the Fredholm alternative holds.
Let us define $R=R(L_{\ov{\mu}})\subseteq C^{r}(\ov{\om}\,)$  to be the range of $L_{\ov{\mu}}$.
Now if $L_{\ov{\mu}}$ is surjective then, by the Fredholm alternative,
we have Ker$(L_{\ov{\mu}})=\emptyset$, which is $(a)$ in the statement of
Theorem \ref{thm1.generic}. Therefore we can assume without loss of generality that
$L_{\ov{\mu}}$ is not surjective, let
$\ov{d}=\mbox{codim}(R)$ be the codimension of $R$.
Since $\ov{L}$ is surjective, by the Fredholm alternative
it is not difficult to see that $\ov{d}\leq 1$, and since $L_{\ov{\mu}}$ is not surjective,
then necessarily $\ov{d}=1$. We will conclude the proof by showing that $(b)$ holds in this case.
Indeed, obviously the kernel must be one dimensional,
 Ker$(\ov{L})=\mbox{\rm span}\{\ov{\phi}\}$, for some $\ov{\phi}\in C^{2,r}_0(\ov{\om}\,)$ which
 satisfies $L_{\ov{\mu}}[\ov{\phi}]=0$. Since  $\ov{L}$ is surjective,
 then $\ov{f}^{'}\ov{\phi}$ must be an element of its range and then there exists
 $\phi\in C^{2,r}_0(\ov{\om}\,)$ which satisfies
$$
L_{\ov{\mu}}[\phi]+\dot{\mu}\ov{f}=\ov{f}^{'}\ov{\phi}.
$$
Multiplying this equation by $\ov{\phi}$ and integrating by parts we find that
$$
\dot{\mu}\ino \ov{f}\,\ov{\phi}=\ino  \ov{f}^{'}\,\ov{\phi}^2,
$$
and since $f^{'}(t)>0$, $\forall\,t$, then we deduce that necessarily $\ino \ov{f}\,\ov{\phi}\neq 0$.
In other words $(b)$ of Theorem \ref{thm1.generic} holds and the proof is concluded.
\finedim

\bigskip
\bigskip

We are ready to present the proof of Theorem \ref{Rab}.\\
{\it Proof of Theorem \ref{Rab}.}\\
It is well known (\cite{CrRab}) that, due to ${\bf (H1)}$,
there exists $\mu_*<+\ii$ such that $\mu\leq \mu_*$ for any solution
of \rife{pm} and in  particular that there exists a continuous simple curve of solutions of
\rife{pm} (the branch of minimal solutions)
for any $\mu<\mu_*$ which emanates from $(\mu,\vm)=(0,0)$, which we denote by $\mathcal{G}(\om)$.
In particular, with the notations of section \ref{sec2}, $\mathcal{G}(\om)$ is characterized by
the fact that the first eigenvalue of the linearized equation, which we denote by $\sg_1(\mu,\vm)$,
satisfies $\sg_1(\mu,\vm)>0$ for any $(\mu,\vm)\in \mathcal{G}(\om)$. In view of
${\bf (H2)}$ and standard elliptic theory we have that $v_*=\left.\vm\right|_{\mu=\mu_*}$ is a classical solution and
$\sg_1(\mu_*,v_*)=0$. By Theorem \ref{thm1.generic} we have that $(b)$ holds for $(\mu_*,v_*)$
and then by Proposition
\ref{pr3.1} we can continue $\mathcal{G}(\om)$ to a
continuous and simple curve without bifurcation points, $[0,s_1+\dt_1)\ni s\mapsto (\mu(s),v(s))$
which locally around any point $s_0>0$ admits a real analytic reparametrization, that is, an injective and continuous map
$\gamma_0:(-1,1)\to (s_0-\eps, s_0+\eps)$, $s=\gamma_0(t)$, such that $\gamma_0(0)=s_0$ and
$(\mu(\gamma_0(t)),v(\gamma_0(t)))$ is real analytic.
Therefore locally this branch has also the structure of a $1$-dimensional
real analytic manifold and we denote it by,
$$
\mathcal{G}^{(s_1+\dt_1)}=\left\{[0,s_1+\dt_1)\ni s\mapsto (\mu(s),v(s))\right\},
$$
which satisfies,
$$
\ov{\mathcal{G}^{(s_1+\dt_1)}}=\left\{[0,s_1+\dt_1]\ni s\mapsto (\mu(s),v(s))\right\},
$$
where, for some $s_1>0$ and $\dt_1>0$, we have:\\
$(A1)_0$ $(\mu(s),v(s))$ is continuous and locally (up to reparametrization) real analytic
for $s\in [0,s_1+\dt_1]$;\\
$(A1)_1$ $v(s)$ is a solution of \rife{pm} with $\mu=\mu(s)$ for any $s\in [0,s_1+\dt_1]$;\\
$(A1)_2$ $\mu(s)=s$ for $s\leq s_1$, $\mu(s_1)=\mu_*$;\\
$(A1)_3$ the inclusion $\left\{(\mu,\vm),\,\mu \in [0,\mu_*] \right\}
\equiv \ov{\mathcal{G}{(\om)}}\subset \mathcal{G}^{(s_1+\dt_1)}$,
holds;\\
$(A1)_4$ $\inf\limits_{[s_1,s_1+\dt_1)}\mu(s)>0$ and
$0<\mu(s)\leq \mu_*,\, \forall\,s\in (0,s_1+\dt_1)$;\\
$(A1)_5$ $0\notin \Sigma(L_{\sscp \mu(s)}),\, \forall\,s\in (0,s_1+\dt_1)\setminus\{s_1\}$,\\
$(A1)_6$ $Ker(L_{\sscp \mu(s_1)})=\,$span$\{\phi_1\}$ and $\ino f(v(s_1))\phi_1\neq 0$,\\
where $\Sigma(L_{\sscp \mu(s)})$ denotes the spectrum of $L_{\sscp \mu(s)}$. Clearly $(A1)_6$
follows from $(b)$ of Theorem \ref{thm1.generic}. Concerning $(A1)_5$ we recall that, by
Proposition \ref{pr3.1}, either $\sg_1(s)$ vanishes identically around $s_1$
or its zero must be isolated. In particular,
since $\sg_1(s)$ is (locally up to a reparametrization) real analytic, its level sets
cannot have accumulation points unless $\sg_1(s)$ is locally constant and consequently
unless it is constant on $[0,s_1+\dt_1)$. However we can rule out this case since,
in view of $(A1)_2$, for $s<s_1$ we have $\sg_1(\mu(s),v(s))>0$
and then no $\sg_k(s)$ can vanish identically, which shows that $(A1)_5$ holds as well. Therefore it is well defined,
$$
s_2:=\sup\left\{t>s_1\,:\inf\limits_{s\in[s_1,t)}\mu(s)>0,\,0\notin \Sigma(L_{\sscp \mu(s)}),\;
\forall\,(\mu(s),v(s))\in \mathcal{G}^{(t)},\,\forall s_1<s<t\right\}.
$$

At this point either $\inf\limits_{s\in[s_1,s_2)}\mu(s)=0$ or $\inf\limits_{s\in[s_1,s_2)}\mu(s)>0$.\\
If $\inf\limits_{s\in[s_1,s_2)}\mu(s)=0$ we set $s_\ii=s_2$,
\beq\label{ii2}
\mathcal{G}^{(s_\ii)}=\left\{[0,s_\ii)\ni s\mapsto (\mu(s),v(s))\right\},
\eeq
and claim that in this case necessarily $\mu(s)\to 0^+$ and $\|v(s)\|_{\ii}\to +\ii$ as $s\to s_\ii$.\\
We first prove that $\mu(s)\to 0^+$ and argue by contradiction, assuming that there exists a
sequence $\{s_n\}\subset (0,s_\ii)$ such that
$s_n\to s_\ii$, $\ainf$ and
$\mu(s_n)\geq \dt>0$ for some $\dt>0$. In view of $(A1)_4$, passing to a subsequence if necessary,
we can assume that $\mu(s_n)\to \ov{\mu}\in [\dt,\mu_*]$.  By ${\bf (H2)}$ and passing to a further
subsequence we would deduce that $v(s_n)\to \ov{v}$, where $(\ov{\mu},\ov{v})$ is a solution of
\rife{pm}. By Theorem \ref{thm1.generic} we see that either $(a)$ or $(b)$ holds and then,
possibly with the aid of Proposition \ref{pr3.1}, we would deduce that locally around
$(\ov{\mu},\ov{v})$
the set of solutions of \rife{pm} is a real analytic parametrization of the form
$(\ov{\mu}(t),\ov{v}(t))$, $t\in (-\eps,\eps)$ for some $\eps>0$
with $(\ov{v}(0),\ov{v}(0))=(\ov{\mu},\ov{v})$.
In particular, for any fixed $\ov{n}$ large enough we can assume without loss of generality that
$(\ov{\mu}(t),\ov{v}(t))$,
$t\in (-\eps,0)$  coincides with $(\mu(s),v(s))$, $s\in (s_{\ov{n}},s_{\ii})$. Now by construction
$\mu(s)>0$ in $[s_1,s_{\ov{n}}]$, and since $\mu(s)$ is continuous we have
$\inf\limits_{s\in [s_1,s_{\ov{n}}]}\mu(s)\geq \ov{\dt}>0$
for some $\ov{\dt}>0$. On the other side, possibly taking a larger $s_{\ov{n}}$, we have
$\inf\limits_{s\in [s_{\ov{n}},s_{\ii})}\mu(s)\geq \frac{\dt}{2}$. In other words we have a contradiction
to $\inf\limits_{s\in [s_1,s_{\ii})}\mu(s)=0$ and the claim is proved.\\
Next we show that $\|v(s)\|_{\ii}\to +\ii$ and argue by contradiction. If this was not the case
we could find a
sequence $\{s_n\}\subset (0,s_\ii)$ such that
$s_n\to s_\ii$, $\ainf$ and $\|v(s_n)\|_{\ii}\leq C$ for some $C>0$. Since we have shown that $\mu(s)\to 0^+$
as $s\to s_{\ii}$, then passing to a subsequence we would deduce that $v(s_{n_k})\to \ov{v}$ where
$\ov{v}$ solves \rife{pm} with $\mu=0$. However by $(A1)_4$ this fact implies that $(\mu,\vm)=(0,0)$
would be a bifurcation point which is clearly impossible, which proves the claim. At this point,
since by definition $\mathcal{R}_{\ii}$ is a closed and connected set, it is not
difficult to see that $\mathcal{R}_{\ii}\equiv \mathcal{G}^{(s_\ii)}$.\\
After a suitable reparametrization we can assume without loss of generality that $s_2=+\ii$ and we conclude that
statement of Theorem \ref{Rab} is true as far as $\inf\limits_{s\in(0,s_2)}\mu(s)=0$. Therefore
we can assume without loss of generality that $\inf\limits_{s\in(0,s_2)}\mu(s)>0$. In this
case, in view of
$(A1)_4$, ${\bf (H2)}$, Theorem \ref{thm1.generic} and Proposition \ref{pr3.1},
it is not difficult to see that $(\mu(s),v(s))$ converges to a solution $(\mu_2,v_2)$ as
$s\to s_2$ and that $0\in \Sigma(L_{\sscp \mu_2})$ and in particular that we can continue the branch
$\mathcal{G}^{(s_2)}$ in a right neighborhood of $s_2$ to a continuous curve which admits
local real analytic reparametrizations. In particular, by arguing as above we see that
$0\notin \Sigma(L_{\sscp \mu(s)})$ for $s\notin \{s_1,s_2\}$ and we can argue by induction
defining, for $k\geq 3$,
$$
s_k:=\sup\left\{t>s_{k-1}\,:\,\inf\limits_{s\in[s_{1},t)}\mu(s)>0,\,0\notin \Sigma(L_{\sscp \mu(s)}),\;
\forall\,(\mu(s),v(s))\in \mathcal{G}^{(t)},\,\forall s_{k-1}<s<t\right\}.
$$
If there exists some $k\geq 3$
such that $\inf\limits_{s\in(0,s_k)}\mu(s)=0$, then as for \rife{ii2} we are done.
Otherwise by using $(A1)_4$, ${\bf (H2)}$, Theorem \ref{thm1.generic} and Proposition \ref{pr3.1},
we can find sequences $s_k$ and $\dt_k>0$ such that, for any $k\in \mathbb{N}$ we have,
$s_{k+1}>s_{k}>\cdots>s_2>s_1$, $s_{k}+\dt_{k}<s_{k+1}$ and,\\

$(Ak)_0$ $(\mu(s),v(s))$ is continuous and simple curve without bifurcation points
(which admits local real analytic
reparamterizations) defined for $s\in [0,s_{k}+\dt_{k}]$;\\
$(Ak)_1$ $v(s)$ is a solution of \rife{pm} with $\mu=\mu(s)$ for any $s\in [0,s_k+\dt_k]$;\\
$(Ak)_2$ $\mu(s)=s$ for $s\leq s_1$, $\mu(s_1)=\mu_*$;\\
$(Ak)_3$ the inclusion $\left\{(\mu(s),v(s)),\,s \in [0,s_{k}] \right\}
\equiv \ov{\mathcal{G}^{(s_k)}{(\om)}}\subset \mathcal{G}^{(s_k+\dt_k)}$,
holds;\\
$(Ak)_4$ $\inf\limits_{[s_1,s_k+\dt_k)}\mu(s)>0$ and
$0<\mu(s)\leq \mu_*,\, \forall\,s\in (0,s_k+\dt_k)$;\\
$(Ak)_5$ $0\notin \Sigma(L_{\sscp \mu(s)}),\, \forall\,s\in (0,s_k+\dt_k)\setminus\{s_1,s_2,\cdots,s_{k}\}$,\\
$(Ak)_6$ $Ker(L_{\sscp \lm(s_k)})=\,$span$\{\phi_k\}$ and $\ino f(v(s_k))\phi_k\neq 0$.\\

Let $s_\ii=\lim\limits_{k\to +\ii}s_k$, we claim that:\\
{\bf Claim:} $\mu(s)\to 0^+$ as $s\to s_\ii$.\\
We argue by contradiction and assume that
along an increasing sequence $\{\widehat{s}_j\}$ such that $\widehat{s}_j\to s_{\ii}$,
it holds $\mu(\widehat{s}_j)\geq \dt>0$ for some $\dt>0$. Clearly we can extract a subsequence
$\{s_{k_j}\}\subset \{s_k\}$ such that $s_{k_j}<\widehat{s}_j\leq s_{k_{j+1}}$.
By $(Ak)_4$ and ${\bf (H2)}$ we can extract an increasing subsequence (which we will not relabel) such that
$(\mu(\widehat{s}_{j}),v(\widehat{s}_{j}))$ converges to a solution $(\widehat{\mu},\widehat{v})$
of \rife{pm} as $j\to +\ii$, where $\dt\leq \widehat{\mu}\leq \mu_*$.

By Theorem \ref{thm1.generic} we can apply either Lemma \ref{lem1.1} or
Proposition \ref{pr3.1} and conclude that locally around $(\widehat{\mu},\widehat{v})$
the set of solutions of \rife{pm} is a real analytic parametrization of the form
$(\widehat{\mu}(t),\widehat{v}(t))$, $t\in (-\eps,\eps)$ for some $\eps>0$
with $(\widehat{\mu}(0),\widehat{v}(0))=(\widehat{\mu},\widehat{v})$.
In particular for $j$ large enough we can assume without loss of generality that
$(\widehat{\mu}(t),\widehat{v}(t))$,
$t\in (-\eps,0)$  coincides with $(\mu(s),v(s))$, $s\in (\widehat{s_{j}},s_{\ii})$.
Let $\{\widehat{\sg}_n\}_{n\in\N}$ be the eigenvalues corresponding to
$(\widehat{\mu},\widehat{v})$ and $\{\widehat{\sg}_n(t)\}_{n\in\N}$ be those corresponding to
$(\widehat{\mu}(t),\widehat{v}(t))$.
On one side, since by construction
$0\in \sg(L_{\sscp {\lm(s_{k_j})}})$ and $s_{k_j}<\widehat{s}_j\leq s_{k_{j+1}}$ for any $j$,
then we have that $0\in \sg(L_{\sscp \widehat{\lm}})$.
Indeed, if this was not the case, then, by Lemma \ref{lem1.1} and since the
eigenvalues are isolated, we would have
that there exists a fixed full neighborhood of $0$ with empty intersection with
$\sg(L_{\sscp {\lm(s_{k_j})}})$ for any $j$ large enough, which is a contradiction since the
number of negative eigenvalues is, locally around each positive solution, uniformly bounded.
As a consequence there exists $n\in \N$ such that $\widehat{\sg}_n=0$.
On the other side, since $\widehat{\sg}_n(t)$ is in particular a continuous function of $t$,
by using once more the fact that the eigenvalues are isolated,
possibly passing to a further subsequence if necessary, we must obviously have
$\widehat{\sg}_n(\widehat{t_j})=0$ for some $\widehat{t_j}\to 0^{-}$ as $j\to +\ii$.
Whence $\widehat{\sg}_n$ must vanish identically in $(-\eps,0]$. In particular the $n$-th eigenvalue
relative to $(\mu(s),v(s))$ must vanish identically for $s\in (\widehat{s_j},s_{\ii})$ and therefore
in $[0,s_{\ii})$.
This is again a contradiction to $(Ak)_2$
since for $s<s_1$ we have $\sg_1(\mu(s),v(s))>0$ and then
no eigenvalue can vanish identically.
Therefore a contradiction arises which shows that $\mu(s)\to 0^+$ as $s\to s_\ii$.\\
At this point, arguing as above, it is not difficult to see that $\|v(s)\|_{\ii}\to +\ii$ as $s\to s_\ii$ and, defining
$$
\mathcal{G}^{(s_\ii)}=\left\{[0,s_\ii)\ni s\mapsto (\mu(s),v(s))\right\},
$$
that $\mathcal{R}_{\ii}\equiv \mathcal{G}^{(s_\ii)}$.
After a suitable reparametrization we can assume without loss of generality that $s_2=+\ii$
which concludes the proof.\finedim


\bigskip
\bigskip






\begin{thebibliography}{99}





\bibitem{BJ} D. Bartolucci, A. Jevnikar, {\em On the global bifurcation diagram of the Gel'fand problem},
Analysis and P.D.E. \textbf{14-8} (2021), 2409-2426.



\bibitem{bdt} B. Buffoni, E.N. Dancer, J.F. Toland, {\em The sub-harmonic bifurcation of Stokes waves},
Arch. Rat. Mech. Anal. {\bf 152}(3) (2000), 24-271.

\bibitem{but} B. Buffoni, J. Toland, {Analytic Theory of Global Bifurcation}, (2003) Princeton Univ. Press.








\bibitem{CrRab} M. G. Crandall, P. H. Rabinowitz, {\em Some Continuation and Variational
Methods for Positive Solutions of Nonlinear Elliptic Eigenvalue Problems},
{Arch. Rat. Mech. An.} {\bf  58} (1975), 207-218.



\bibitem{Dan} N. Dancer, {\em Real analyticity and non-degeneracy},
Math. Ann. {\bf 325} (2003), 369-392.

\bibitem{Dan3} N. Dancer,
{\em Global structure of the solutions of non-linear real analytic eigenvalue problems},
Proc. Lond. Math. Soc. (3) {\bf 27} (1973), 747-765.

\bibitem{Dan1} E.N. Dancer, {\em On the structure of solutions of an equation in catalysis theory when
a parameter is large}, J. Diff Eq {\bf 37}(3) (1980),  404-437.



\bibitem{DLN} D.G. de Figuereido, P.-L. Lions, R.D. Nussbaum
{\em A priori estimates and existence of positive solutions of semilinear elliptic equations},
J. Math. Pures Appl. {\bf 61}(1) (1982), 41-63.







\bibitem{JL} D.D. Joseph, T.S. Lundgren, {\em
Quasilinear Dirichlet problems driven by positive sources},
Arch. Ration. Mech. Anal. {\bf 49} (1972/73) 241-269.


\bibitem{HK}
M. Holzmann, H. Kielh\"ofer, {\em Uniqueness of global positive solution branches
of nonlinear elliptic problems}, Math. Ann. {\bf  300}, 221-241 (1994).


\bibitem{H} D. Henry, {\sl Perturbation of the boundary in boundary-value problems of partial differential equations}, London Math. Soc.
L.N.S. {\bf (318)} Cambridge University Press, Cambridge, (2005).



\bibitem{Ka} T. Kato, {\sl Perturbation Theory of Linear Operators}, Springer-Verlag (1966).

\bibitem{Kor}  P. Korman, {\em Global Solution Curves for Semilinear Elliptic Equations}, World Scientific (2012).


\bibitem{Kor1} P. Korman, {\em
Global solution curves for self-similar equations}, J. Diff. Eq. {\bf 257} (2014) 2543-2564.





\bibitem{Mich}
A. M. Micheletti, {\em Perturbazione dello spetro di un operatore ellitico di tipo variazionale,
in relazione ad una variazione del campo}, Annali Mat. Pura App. {\bf 97} (1973), 267-281.


\bibitem{N}  S. Nakane {\em A bifurcation phenomenon for a Dirichlet problem with an exponential nonlinearity},
J. Math. An. App. 161 (1991), 227-240.

\bibitem{NS}  K. Nagasaki, T. Suzuki, {\em Radial and nonradial solutions for the nonlinear eigenvalue problem
$\Delta u + \lm e^u = 0$ on annuli in $\mathbb{R}^2$}, J. Differential Equations 87 (1990), 144-168.

\bibitem{NS90} K. Nagasaki, T. Suzuki, {\em Asymptotic analysis for two-dimensional elliptic eigen-
values problems with exponentially dominated nonlinearities}, Asymptotic Analysis,
{\bf 3} (1990), 173-188.


\bibitem{PW} M Plum, C. Wieners,  {\em New solutions of the Gelfand problem},
J. Math. An. App. 269 (2002), 588-606.




\bibitem{Rab} P.H. Rabinowitz, {\em Some global results for nonlinear eigenvalue problems}, J. Funct. Anal.
{\bf 7} (1981), 487-513.

\bibitem{Rynne} B.P. Rynne {\em The structure of Rabinowitz' global bifurcating
continua for generic quasilinear elliptic equations} Noniln. An. T.M.A. {\bf 32} (1998) 167-181.


\bibitem{ST} J.C. Saut, R. Temam, {\em Generic properties of nonlinear boundary-value problems},
Comm. P.D.E. (1979), 293-319.

\bibitem{Sm} S. Smale, {\em An infinite-dimensional version of Sard's theorem}, Amer. J. Math. {\bf 87} (1969), 861-866.


\bibitem{Uhl} K. Uhlenbeck, {\em Generic properties of eigenfunctions}, Amer. J. Math. {\bf 98} (1976),
1059-1078.


\end{thebibliography}
\end{document}